\documentclass[letterpaper,10pt]{article}

\usepackage[margin=1in]{geometry}
\usepackage[utf8]{inputenc}
\usepackage{amsmath}
\usepackage{mathtools}
\usepackage{amsthm,amssymb}
\usepackage{amsfonts}
\usepackage[numbers,square]{natbib}
\usepackage{comment}
\usepackage{authblk}

\newtheorem{theorem}{Theorem}

\numberwithin{equation}{section}

\newcommand{\id}{ {1\!\!1 } }

\DeclareMathOperator{\Log}{Log}
\DeclareMathOperator{\GL}{GL}
\DeclareMathOperator{\SO}{SO}
\DeclareMathOperator{\Oo}{O}
\renewcommand{\skew}{\mathop{\rm skew}}
\DeclareMathOperator{\sym}{sym}
\DeclareMathOperator{\dev}{dev}

\DeclareMathOperator{\so}{\mathfrak{so}}
\DeclareMathOperator{\gl}{\mathfrak{gl}}
\DeclareMathOperator{\polar}{polar} 
\DeclareMathOperator{\GLp}{GL^+}
\DeclareMathOperator{\dist}{dist}
\DeclareMathOperator{\Psym}{PSym}
\DeclareMathOperator{\tr}{tr}

\newcommand{\R}{\mathbb{R}}
\newcommand{\Rnn}{\R^{n\times n}}
\newcommand{\ds}{\mathrm{ds}}
\newcommand{\tel}[1]{\frac{1}{#1}}
\newcommand{\half}{\tel{2}}

\newcommand{\norm}[1]{|\mkern-2mu|#1|\mkern-2mu|}
\newcommand{\innerproduct}[1]{\left\langle #1 \right\rangle}

\newcommand{\len}{L}
\newcommand{\dg}{\dist_{\rm geod}}

\newcommand{\fullnorm}[1]{{\norm{#1}}_{\mu, \mu_c, \kappa}}
\newcommand{\fullproduct}[1]{\innerproduct{#1}_{\mu, \mu_c, \kappa}}
\newcommand{\adset}{\mathcal{A}}

\newcommand{\generalproduct}[1]{\innerproduct{#1}}
\newcommand{\omegatfrac}{\tfrac{\mu_c}{\mu}}
\newcommand{\eucproduct}[1]{\innerproduct{#1}_{n\times n}}

\newcommand{\smallnorm}[1]{\norm{\mkern-2mu #1 \mkern-1mu}}
\newcommand{\smallfullnorm}[1]{\fullnorm{\!#1\mkern-1mu}}

\DeclareMathOperator{\Curl}{Curl}

\begin{document}

\title{A Riemannian approach to strain measures in nonlinear elasticity}
\author[1]{Patrizio Neff$\mkern1mu$\footnote{To whom correspondence should be addressed. e-mail: \texttt{patrizio.neff@uni-due.de}}}
\author[2]{Bernhard Eidel}
\author[1]{Frank Osterbrink}
\author[1]{Robert Martin}
\affil[1]{Chair for Nonlinear Analysis and Modelling, University of Duisburg-Essen, Thea-Leymann-Str. 9, 45127 Essen}
\affil[2]{Institute of Mechanics, University of Duisburg-Essen, Universitätsstraße 15, 45141 Essen}
\date{\today}
\maketitle

\begin{abstract}
	The isotropic Hencky strain energy appears naturally as a distance measure of the deformation gradient to the set $\SO(n)$ of rigid rotations in the canonical left-invariant Riemannian metric on the general linear group $\GL(n)$. Objectivity requires the Riemannian metric to be left-$\GL(n)$-invariant, isotropy requires the Riemannian metric to be right-$\Oo(n)$-invariant. The latter two conditions are satisfied for a three-parameter family of Riemannian metrics on the tangent space of $\GL(n)$. Surprisingly, the final result is basically independent of the chosen parameters.

In deriving the result, geodesics on $\GL(n)$ have to be parametrized and a novel minimization problem, involving the matrix logarithm for non-symmetric arguments, has to be solved.
\end{abstract}

\section{Introduction}
For the deformation gradient $F = \nabla \varphi \in \GLp(n)$ let $U=\sqrt{F^TF}$ be the symmetric right Biot-stretch tensor.
We show that the isotropic Hencky strain energy, defined on the logarithmic strain tensor $\log U$ by
\begin{align}
	W(F) = \mu \, \smallnorm{\dev\log U}^2+\frac{\kappa}{2}[\tr(\log U)]^2 = \mu \, \smallnorm{\dev\log U}^2+\frac{\kappa}{2}(\log \det F)^2\,,
\end{align}
measures the geodesic distance of $F$ to the group of rotations $\SO(n)$ where $\GL(n)$ is viewed as a Riemannian manifold endowed with a left-invariant metric which is also right $\Oo(n)$-invariant (isotropic), and where the coefficients \mbox{$\mu,\, \kappa > 0$} correspond to the shear modulus and the bulk modulus, respectively. Thus we provide yet another characterization of the polar decomposition $F = R \, U\, ,\, R\in\SO(n),\, U\in \Psym(n)$, since $U$ also provides the minimal euclidean distance to $\SO(n)$, i.e.,
\begin{align}
	\dist_{\rm euclid}^2(F,\SO(n)) := \!\!\!\underset{Q\in\SO(n)}{\min} \!\dist_{\rm euclid}^2(F,Q) = \!\!\!\underset{Q\in\SO(n)}{\min} \!\norm{F-Q}^2 = \norm{F-R}^2 = \norm{U-\id}^2\,, \label{eq:introductionEuclideanDistanceToSO}
\end{align}
where the euclidean distance $\,{\rm{dist}}_{\rm{euclid}}^2(X,Y):=\norm{X-Y}^2$ is the length of the line segment joining $X$ and $Y$ in $\R^{n^2}$, $\id\in\GLp(n)$ is the identity and $\norm{X} = \sqrt{\tr(X^T X)}$ denotes the Frobenius matrix norm here and henceforth.
For both the euclidean and the geodesic distance, the orthogonal factor $R = \polar (F)$ in the polar decomposition of $F$ is the nearest rotation to $F$.

\section{Strain measures in linear and nonlinear elasticity}
We consider an elastic body which in a reference configuration occupies the bounded domain $\Omega\subset\R^3$. Deformations of the body are prescribed by mappings $\varphi:\Omega\to \R^3$, where $\varphi(x)$ denotes the deformed position of the material point $x\in\Omega$. Central to elasticity theory is the notion of strain, which is a measure of deformation such that vanishing strain implies that the body $\Omega$ has been moved rigidly in space.
Various such measures exist, e.g. the \emph{Green strain} $\half (U^2-\id)$, the \emph{generalized Green strain} $\tel{m}(U^m-\id)$, where $m$ is a nonzero integer, and the \emph{Hencky (or logarithmic) strain} $\log U$.

In linearized elasticity, one considers $\varphi(x)=x+u(x)$, where $u:\Omega\subset\R^3\to \R^3$ is the displacement. The classical linearized strain measure is $\varepsilon=\sym\nabla u$. It appears through a matrix-nearness problem in the euclidean distance
\begin{align}
    \dist_{\rm euclid}^2(\nabla u,\so(3)):=\min_{W\in\so(3)} \norm{\nabla u-W}^2=\norm{\sym \nabla u}^2\,,
\end{align}
where $\so(3)$ denotes the set of all skew symmetric matrices in $\R^{3 \times 3}$. Indeed, $\sym\nabla u$ qualifies as a linearized strain measure: if $\dist_{\rm euclid}^2(\nabla u,\so(3))=0$ then $u(x)=\widehat{W}.x+\widehat{b}$ is a linearized rigid movement. This is the case since 
\begin{align}
      \dist_{\rm euclid}^2(\nabla u(x),\so(3))=0 \quad \Rightarrow \quad \nabla u(x)=W(x)\in\so(3)
      \end{align}
and $0=\Curl \nabla u(x)=\Curl W(x)$ implies that $W(x)$ is constant, see \cite{Neff_curl06}.
In nonlinear elasticity theory one assumes that $\nabla \varphi\in\GL^+(3)$ (no self-interpenetration of matter) and may consider the matrix nearness problem
\begin{align}
    {\rm{dist}}_{\rm{euclid}}^2(\nabla \varphi,\SO(3)):=\min_{Q\in\SO(3)} \norm{\nabla \varphi-Q}^2=\min_{Q\in\SO(3)}\norm{Q^T\nabla\varphi-\id}^2 = \norm{\sqrt{\nabla \varphi^T \nabla\varphi}-\id}^2\, , \label{dist_euclid}
\end{align}
where the last equality is due to \eqref{eq:introductionEuclideanDistanceToSO}. Indeed, the Biot strain tensor $\sqrt{\nabla \varphi^T \nabla\varphi}-\id$ qualifies as a nonlinear strain measure: if $ {\rm{dist}}_{\rm{euclid}}^2(\nabla\varphi,\SO(3) ) =0$ then $\varphi(x)=\widehat{Q}.x+\widehat{b}$ is a rigid movement. This is the case since 
\begin{align}
     {\rm{dist}}_{\rm{euclid}}^2(\nabla\varphi,\SO(3) )=0 \quad \Rightarrow \quad \nabla\varphi(x)=Q(x)\in\SO(3)
      \end{align}
and $0=\Curl \nabla \varphi(x)=\Curl Q(x)$ implies that $Q(x)$ is constant, see \cite{Neff_curl06}.

In geometrically nonlinear, physically linear isotropic elasticity the formulation of a boundary value problem of place may now be based on minimizing the quadratic Biot strain energy
\begin{align}
	\mathcal{E}(\varphi) = \int_\Omega \mu\, \norm{\dev [\sqrt{\nabla\varphi^T\nabla\varphi}-\id]}^2+\frac{\kappa}{2}\left( \tr{\sqrt{\nabla\varphi^T\nabla\varphi}-\id}\right)^2\,{\rm dx}\,, \quad \left.\varphi\right|_{\Gamma_D}=\varphi_0\,,
\end{align}
where $\mu, \kappa > 0$ are the shear modulus and bulk modulus, respectively.

However, since the Euclidean distance in \eqref{dist_euclid} is an arbitrary choice, novel approaches in nonlinear elasticity theory aim at putting more geometry (i.e. respecting the group structure of the deformation mappings) into the description of the strain a material endures.  In our context, it is now natural to consider a strain measure induced by the geodesic distances stemming from choices for the Riemannian structure respecting also the algebraic group structure of $\GLp(n)$, which we introduce next.

\section{Left invariant Riemannian metrics on $\GL(n)$}

Viewing $\GL(n)$ as a Riemannian manifold endowed with a left invariant metric
\begin{align}
	g_H: T_H \GL(n) \times T_H \GL(n) \to \R: \:\: g_H(X,Y) = \generalproduct{H^{-1}X,\: H^{-1}Y}\,,\quad H\in\GL(n)\,,
\end{align}
for a suitable inner product $\generalproduct{\cdot,\cdot}$ on the tangent space $T_\id \GL(n) = \gl(n) = \Rnn$ at the identity $\id$, the distance between $F,P\in\GLp(n)$ can be measured along sufficiently smooth curves. We denote by
\begin{align}
	\adset = \{\gamma\in C^0([0,1];\GLp(n))\: | \:\:\gamma\: \text{ piecewise differentiable, }\, \gamma(0)=F,\, \gamma(1)=P\}
\end{align}
the admissible set of curves connecting $F$ and $P$, and by
\begin{align}
	\len(\gamma) = \smash{\int\nolimits_0^1} \sqrt{g_{\gamma(s)} (\dot\gamma(s), \dot\gamma(s))}\:\ds
\end{align}
the length of $\gamma\in\adset$. Then the geodesic distance
\begin{align}
	\dg(F,P) = \underset{\gamma\in\adset}{\inf} \len(\gamma)
\end{align}
defines a metric on $\GLp(n)$. While it is generally difficult to explicitly compute this distance or to find length minimizing curves, it can be shown \cite{Mielke02,Neff_Martin_2013} that if the Riemannian metric is defined by an inner product of the form
\begin{align}
	&\generalproduct{X,Y} = \fullproduct{X,Y} := \mu\eucproduct{\dev\sym X, \dev\sym Y} + \mu_c\eucproduct{\skew X, \skew Y} + \tfrac{\kappa}{2}\tr X \tr Y\,, \nonumber \\
	&\fullnorm{X}^2 := \fullproduct{X,X} = \mu\,\smallnorm{\dev\sym X}^2 + \mu_c\,\smallnorm{\skew X}^2 + \tfrac{\kappa}{2}\,[\tr X]^2\,,\qquad \mu,\,\mu_c,\,\kappa > 0\,, \label{eq:isotropicProduct} \\
	&\dev X := X - \tfrac{1}{n}\tr X \cdot \id\,,\quad \mu_c \text{ denoting the \emph{spin modulus}}\,, \nonumber
\end{align}
which is the case if and only if the metric $g$ is right invariant under $\Oo(n)$ \cite{Bryant2013}, then every geodesic $\gamma$ connecting $F$ and $P$ is of the form
\begin{align}
	\gamma(t) = F\, \exp(t(\sym\xi - \omegatfrac\skew\xi))\: \exp(t(1 + \omegatfrac)\skew\xi)
\end{align}
for some $\xi\in\gl(n)$, where $\exp: \gl(n) \to \GLp(n)$ denotes the matrix exponential, $\sym \xi = \half(\xi + \xi^T)$ the symmetric part and $\skew \xi = \half(\xi - \xi^T)$ the skew symmetric part of $\xi$.

Now, according to the classical Hopf-Rinow theorem of differential geometry, there exists a length minimizing geodesic in $\adset$ for all $F,P\in\GLp(n)$. To obtain such a minimizer $\gamma$ (and thus the distance $\dg(F,P) = \len(\gamma)$), it therefore remains to find $\xi\in\gl(n)$ with
\begin{align}
	P = \gamma(1) = F\, \exp(\sym\xi - \omegatfrac\skew\xi)\: \exp((1 + \omegatfrac)\skew\xi)\,.\label{eq:geodesicendpoint}
\end{align}
The existence of such a $\xi$ is clear from the above.

\section{The geodesic distance to $\SO(n)$}

Although no closed form solution to \eqref{eq:geodesicendpoint} is known, the equation can be used to obtain a lower bound\footnote{We denote by $\log$ the principal matrix logarithm, while the expression $\Log$ is used to indicate that the infimum is taken over the whole inverse image under $\exp$, i.e. $\underset{Q\in\SO(n)}{\min} \smallfullnorm{\Log(QF)}^2 = \min\{\fullnorm{\xi}^2:\: \xi\in\gl(n),\: \exp(\xi)=QF\}$.}
\begin{align}
	\dg^2(F,\SO(n)) = \underset{Q\in\SO(n)}{\min} \! \dg^2(F,Q) \geq \underset{Q\in\SO(n)}{\min} \smallfullnorm{\Log(QF)}^2 \label{eq:lowerBound}
\end{align}
for the distance of $F\in\GLp(n)$ to $\SO(n)$, as well as a simple upper bound
\begin{align}
	\dg^2(F,\SO(n)) &\leq \dg^2(F,\, \polar(F)) \nonumber \\
	&\leq \smallfullnorm{\log(\polar(F)^T F)}^2 \:=\:\: \mu\smallnorm{\dev\log(U)}^2 + \frac{\kappa}{2}[\tr(\log U)]^2\,, \label{eq:upperBound}
\end{align}
where $F = R\,U$, $R=\polar(F)\in\SO(n)$, $U=\sqrt{F^T F}\in\Psym(n)$ denotes the polar decomposition of $F$.
Finally, we can use an extension of a recent optimality result proved by Neff et al. \cite{Neff_Nakatsukasa_Fischle_2013}:

\begin{theorem}
	\label{theorem:infLog}
	Let $\norm{\,.\,}$ be the Frobenius matrix norm on $\gl(n)$, $F\in\GLp(n)$. Then the minimum
	\begin{align}
		\underset{Q\in\SO(n)}{\min} \smallnorm{\Log(Q \cdot F)}^2 = \smallnorm{\log(\polar(F)^T F)}^2 =  \smallnorm{\log(\sqrt{F^T F})}^2 = \smallnorm{\log(U)}^2\label{eq:infLog}
	\end{align}
	is uniquely attained at $Q=\polar(F)^T$.
\end{theorem}

A consequence of Theorem \ref{theorem:infLog}, combined with \eqref{eq:isotropicProduct}, \eqref{eq:lowerBound} and \eqref{eq:upperBound}, yields our main result \cite{Neff_Eidel_Osterbrink_2013}:

\begin{theorem}
	\label{theorem:mainResult}
	Let $g$ be a left invariant Riemannian metric on $\GL(n)$ that is also right invariant under $\Oo(n)$, and let $F\in\GLp(n)$. Then:
	\begin{align}
		\dg^2(F,\SO(n)) = \dg^2(F,\,\polar(F)) = \mu\smallnorm{\dev\log(U)}^2 + \frac{\kappa}{2}[\tr(\log U)]^2\,.
	\end{align}
	Thus the geodesic distance of the deformation gradient $F$ to $\SO(n)$ is the isotropic Hencky strain energy of $F$. In particular, the result is independent of the spin modulus $\mu_c>0$.
\end{theorem}

Furthermore, for $\mu_c = 0$ (in which case $\dg$ defines only a pseudometric on $\GLp(n)$), Theorem \ref{theorem:mainResult} still holds.

\end{document}